\documentclass{article}

\begin{document}

%+Title
\title{Twistor construction of some multivalued harmonic functions on ${\bf R}^{3}$}
\author{Simon Donaldson}
\date{\today}
\maketitle
%-Title

%+Abstract
%\begin{abstract}
%    There is abstract text that you should replace with your own. 
%\end{abstract}
%-Abstract

%+Contents
%\tableofcontents
%-Contents

\section{Introduction}
\newcommand{\bR}{{\bf R}}
\newcommand{\bC}{{\bf C}}
\newcommand{\bT}{{\bf T}}
\newcommand{\bZ}{{\bf Z}}
\newcommand{\cO}{{\cal O}}
\newcommand{\cF}{{\cal F}}
\newcommand{\cU}{{\cal U}}
\newcommand{\tphi}{\widetilde{\phi}}
\newcommand{\uphi}{\underline{\phi}}
\newcommand{\oA}{\overline{A}}
\newtheorem{thm}{Theorem}
\newtheorem{lem}{Lemma}
\newtheorem{prop}{Proposition}
Over the past decade, multivalued solutions of elliptic partial differential equations have become important in a number of areas of differential geometry. These include questions in gauge theory---starting with the work of Taubes \cite{kn:T}---in the geometry of calibrated submanifolds, as in the work of He \cite{kn:He} and in the area of exceptional holonomy. Here we will only consider harmonic functions. A multivalued harmonic function on a Riemannian manifold $M$ with branch set a codimension $2$ submanifold $\Gamma\subset M$ is a harmonic section of a  flat real line bundle over $M\setminus \Gamma$ with holonomy $-1$ around  small loops linking $\Gamma$. Alternatively, it can be viewed as a function on a double branched cover, anti-invariant under the covering involution.

 If we consider a fixed $\Gamma$ and allow functions which are $O(r^{1/2})$ near $\Gamma$, where $r$ is the distance to $\Gamma$, there is a relatively straightforward theory. From one point of view,  the lift of the Riemannian metric on $M$ to the double cover is a singular metric but the singularity can be handled by established elliptic theory. For the case $M=\bR^{n}$, which is our concern in this paper, it was shown by  by Sun \cite{kn:Sun} and Haydys,Mazzeo, Takahashi \cite{kn:HMT} that for any fixed $\Gamma$ and harmonic polynomial $F$ there is a multivalued harmonic function with branch set $\Gamma$ asymptotic to $F$ at infinity,   But in the applications indicated above one needs functions in a more restricted class, with locally bounded derivative---they will be $O(r^{3/2})$ near $\Gamma$. In  existence questions it is not then usually possible to fix $\Gamma$  from the outset---the submanifold has to be found as part of the solution. The problem is nonlinear and there are rather few general existence results.

 An important advance was made in recent work of Dashen Yan \cite{kn:Yan} who found  explicit solutions on $\bR^{n}$ in all dimensions $n$,  using an ingenious choice of coordinate system. The purpose of this note is to give another construction of a family of explicit solutions in the  case of dimension $n=3$ using a classical integral representation, related to twistor theory. While we restrict attention here to the particular family of solutions, described in Theorem 1, one can hope that similar techniques can be applied to construct other solutions---for various PDE---in dimensions $3$ and $4$.

Take standard Euclidean coordinates $(X,Y,Z)$ on $\bR^{3}$. For $a,b>0$ let
$\Gamma$ be the ellipse
$$   \frac{X^{2}}{a^{2}} +  \frac{Y^{2}}{b^{2}}=1 $$
in the plane $Z=0$. Set
$$  W_{0}= \{ (X,Y,0): X^{2}/a^{2}+Y^{2}/b^{2}<1 \}, $$
and $$   W_{\infty}= \{ (X,Y,0): X^{2}/a^{2}+Y^{2}/b^{2}>1 \}.$$
A multivalued function with branch set $\Gamma$ can be described by a function $\phi$ on $\bR^{3}\setminus \overline{W_{\infty}}$ satisfying matching conditions on the wall $W_{\infty}$. Write $\phi_{Z}$ for the partial derivative $\frac{\partial \phi}{\partial Z}$. We require that $\phi(p)$ and $\phi_{Z}(p)$ have well-defined limits as the point $p$ approaches $W_{\infty}$ from the two half-spaces. Denote these limits (which are functions on $W_{\infty}$)  by $\phi^{+},\phi^{-}$ and $\phi_{Z}^{+}, \phi_{Z}^{-}$.  The matching conditions are
\begin{equation}  \phi^{-}=-\phi^{+} \ \ \ \ , \ \ \  \  \phi_{Z}^{-}=-\phi_{Z}^{+}. \end{equation}
The same multivalued section has another description as a function $\tphi$ on $\bR^{3}\setminus \overline{W_{0}}$ satisfying the similar matching conditions on $W_{0}$. The two descriptions are related by
$$  \tphi(X,Y,Z)= {\rm sign}(Z)\  \phi(X,Y,Z). $$

\begin{thm} For any $\epsilon$ with $-1<\epsilon<1$ set $a=\sqrt{1+\epsilon}, b=\sqrt{1-\epsilon}$. There is a harmonic function $\phi$ on $\bR^{3}\setminus \overline{W_{\infty}}$ satisfying the matching conditions (1) and such that
\begin{enumerate}
 \item $\nabla \phi$ is bounded near $\Gamma$;
 \item  At large distances $$\tphi(X,Y,Z)= -\lambda X^{2}- \mu Y^{2} + \nu Z^{2} +O(R), $$ 
where $R=\sqrt{X^{2}+Y^{2}+Z^{2}}$ and
$$  \lambda=(1/4) \int_{0}^{2\pi} \frac{\cos^{2} \theta}{(1+\epsilon  \cos 2\theta)^{3/2}}d\theta\ \ \ \ \ \ \mu=(1/4) \int_{0}^{2\pi} \frac{\sin^{2} \theta}{(1+\epsilon \cos 2\theta)^{3/2}}d\theta $$
$$  \ \ \  \nu= (1/4) \int_{0}^{2\pi} \frac{1}{(1+\epsilon \cos 2\theta)^{3/2}}d\theta $$
\end{enumerate}
\end{thm}

{\bf Remark.} \ Our solutions $\phi$ will be odd in the $Z$-variable ({\it i.e.} $\phi(X,Y,-Z)=-\phi(X,Y,Z)$). For such functions the matching conditions (1) become $\phi_{Z}^{\pm}=0$ on $W_{\infty}$. The problem can be viewed as a half-space problem with mixed Dirichlet and Neumann boundary conditions: find a harmonic function $\phi$ on the half-space $\{Z\geq 0\}$ with $\phi=0$ on $W_{0}$, $\phi_{Z}=0$ on $W_{\infty}$ and with $\nabla \phi$ bounded near $\Gamma$. 
\section{Background}

There is a classical integral representation formula for the general local solution of the Laplace equation in $\bR^{3}$:
\begin{equation}
\phi(X,Y,Z)= \frac{1}{2\pi} \int_{0}^{2\pi} f(Z+ i(X \cos \theta +Y \sin\theta),\theta) d\theta, \end{equation}
where $f$ is an analytic function of two variables, defined in a suitable region. See \cite{kn:WW}, page 390.
Changing notation , write $w=e^{i\theta}$ and $A= \frac{1}{2}(Y+iX)$. Then the integral (2) is
\begin{equation}\phi(A,Z)= \frac{1}{2\pi i} \int_{S^{1}} F(Z+ A w -\oA w^{-1}, w) \ \frac{dw}{w} \end{equation} where the integration contour $S^{1}$ is the unit circle and $F$ is analytic in a suitable region.

This formula has a more modern and geometric interpretation through twistor theory---more precisely, the three dimensional version of twistor theory developed by Hitchin \cite{kn:H}. Let $\bT$ be the set of oriented lines in $\bR^{3}$. There is an obvious identification of $\bT$ with the total space of the tangent bundle of $S^{2}$ and so $\bT$ is a $2$-dimensional complex manifold. Reversing the orientation of a line defines an anti-holomorphic involution $\sigma:\bT\rightarrow \bT$. The holomorphic sections of $\bT$ ({\it i.e.} holomorphic vector fields on $S^{2}$) form a $3$-dimensional vector space $\bC^{3}$. For each $p\in \bC^{3}$ the corresponding section defines  a holomorphic curve $L_{p}\subset \bT$. Our original Euclidean $3$-space $\bR^{3}$ is recovered as the set of real points in $\bC^{3}$, corresponding to the \lq\lq real''curves $L_{p}$ which are preserved by $\sigma$. For $p\in \bR^{3}$ the curve $L_{p}\subset \bT$ is the set of lines through $p$. 

The twistor description of harmonic functions goes as follows. Let $\cO(-2)$ be the holomorphic line bundle over $\bT$ which is the lift of the canonical bundle of $S^{2}$ by the bundle projection map $\pi:\bT\rightarrow S^{2}$. For any $p$ the restriction of $\cO(-2)$ to the curve $L_{p}$ is naturally identified with the canonical bundle of $L_{p}$. Let $\cU$ be an open set in $\bT$ and $\cF$ be a class in the sheaf cohomology group $H^{1}(\cU, \cO(-2))$. For any curve $L_{p}$ which lies inside $\cU$ we set
\begin{equation}  \phi(p)= {\rm ev}\ \cF\vert_{L_{p}}, \end{equation}
where
$$  {\rm ev}: H^{1}(L_{p}, \cO(-2))\rightarrow \bC$$ is the natural map induced by the identification of $\cO(-2)$ with the canonical bundle. The main statement is that this defines a harmonic function and all harmonic functions arise in this way, at least locally in $\bR^{3}$.

To recover the classical integral formula, recall that if a class in $H^{1}(S^{2}. \cO(-2))$ is represented by a Cech cocycle in the standard way---a holomorphic $1$-form $\alpha$ on a neighbourhood of the unit circle---then the map ${\rm ev}$ is given by contour integration
\begin{equation}  {\rm ev}\ ([\alpha])= \int_{S^{1}} \alpha. \end{equation}
 Let $w$ be a standard complex co-ordinate on $S^{2}$ with the antipodal points $w=0,\infty$ corresponding to the $Z$-axis. The holomorphic vector field ${\bf v}= w\frac{\partial}{\partial w}$ defines a trivialisation of the bundle $\bT$ away from these two points and we get complex co-ordinates $(w,u)$ on an open set in $\bT$ (with $(w,u)$ corresponding to $u \ {\bf v}(w)\in \bT$).  In these coordinates the antiholomorphic involution is given by $\sigma(w,u)= (-\overline{w}^{-1}, \overline{u})$.

The general holomorphic vector field on $S^{2}$ is $(\alpha w + \beta +\gamma
w^{-1}){\bf v}$ for $(\alpha,\beta,\gamma)\in \bC^{3}$. This defines a real curve if and only if $\beta$ is real and $\gamma=-\overline{\alpha}$. Under our identification of $\bR^{3}$ with the real curves we write $\alpha=A, \beta=Z, \gamma=-\oA$. Thus the curve corresponding to parameters $A,Z$ is given by the equation \begin{equation} u= A w + Z -\oA w^{-1} . \end{equation}
In our coordinates the $1$-form $dw/w$ defines a trivialisation of $\cO(-2)$. If a class in $H^{1}(\cU, \cO(-2))$ is represented by a Cech cocycle
$F(w,u) dw/w$ (with respect to  a suitable cover of $\cU$ by two open subsets) then the two formulae (4),(5) reduce to the classical representation (3) for the corresponding harmonic function.

\section{Proof of Theorem 1}

Fix $\epsilon$ in the interval $(-1,1)$ and define the function
$$Q(w)= 1+ \frac{\epsilon}{2} (w^{2}+ w^{-2}). $$

Define the function $\tan^{-1}$ (also written $\arctan$) by taking two cuts on the imaginary axis in $\bC$ one from $i$ to $+\infty\  i$ and one from $-i$ to $-\infty\  i$. This is an odd function on the cut plane. Now define
\begin{equation} F(w,u)=  Q(w)^{-3/2} (Q(w)+u^{2})\tan^{-1}(u/\sqrt{Q(w)}) + \kappa_{\epsilon} u\end{equation}
where $\kappa_{\epsilon}$ is a constant to be fixed below. We want to define our harmonic function $\phi$ using the formula (3) with this choice of $F$. To do this we need to check first that $F$ is well-defined for the relevant values of $(w,u)$. There is no difficulty with the fractional powers of $Q$ because $\vert Q(w)\vert <1$ for $w$ in the unit circle. To deal with the inverse tangent we need:
\begin{lem}
Suppose $p\in \bR^{3}$ and, with notation as above, that there is a point $w$ on the unit circle $C$ such that 
$$  A w + Z -\oA w^{-1}= t\  i \ \sqrt{Q(w)} $$
for $t$ real with  $\vert t\vert\geq 1$. Then $p$ lies in $\overline{W_{\infty}}\subset \bR^{3}$. 
\end{lem}
This is completely elementary. For $w$ in the  unit circle  $w^{-1}= \overline{w}$. Squaring the equation and re-arranging we get
$$   2 {\rm Re}\ \left( (A^{2}+\frac{\epsilon t^{2}}{2}) w^{2}\right) + t^{2}+ Z^{2} - 2\vert A\vert^{2}= 4 i Z\  {\rm Im}\ \left( A w\right). $$
The left hand side is real and the right hand is pure imaginary so both must vanish. Suppose ${\rm Im}\ (Aw)=0$. Then $Aw$ is real and $A^{2} w^{2}= \vert A\vert^{2}$. So we get
  $$   {\rm Re}\ ( \epsilon t^{2} w^{2}) + t^{2}+ Z^{2}=0. $$
  But this is impossible with our range of $\epsilon$ since
  $\vert \epsilon t^{2} w^{2}\vert < t^{2}$. Thus we deduce that $Z=0$ and
  $$     {\rm Re}\ \left( (A^{2}+\frac{\epsilon t^{2}}{2}) w^{2}\right)
 = \vert A\vert^{2}-t^{2}/2. $$
 This implies that
 $$  \vert (A^{2}+\frac{\epsilon t^{2}}{2}) w^{2}\vert^{2}\geq \left( \vert A\vert^{2}- t^{2}/2\right)^{2}, $$
which gives (using $\vert w\vert=1$):
   $$  \vert A\vert^{4} + {\rm Re}\ (\epsilon t^{2} A^{2}) +\epsilon^{2}t^{4}/4 \geq \vert A\vert^{4} -\vert A\vert^{2} t^{2} + t^{4}/4 . $$
Re-arranging, and dividing by $t^{2}$, we get
$$\vert A\vert^{2} + \epsilon {\rm Re}\ ( A^{2}) \geq \frac{t^{2}}{4} (1-\epsilon^{2}). $$
Recall that $A= \frac{1}{2} (Y-iX)$ so this is:
$$  \frac{X^{2}+Y^{2}}{4}+ \frac{\epsilon}{4}(Y^{2}-X^{2}) \geq \frac{t^{2}}{4} (1-\epsilon^{2}),$$ which gives
$$  (1-\epsilon)X^{2} + (1+\epsilon )Y^{2} \geq t^{2}(1-\epsilon^{2})\geq (1-\epsilon^{2}). $$
So $X^{2}/(1+\epsilon) + Y^{2}/(1-\epsilon) \geq 1$ and the point $p$ lies in $\overline{W_{\infty}}$ as claimed.

We now have a harmonic function $\phi$ on $\bR^{3}\setminus \overline{W_{\infty}}$. By changing the integration variable to $-w$ or $w^{-1}$ one sees that the function has the symmetries:
\begin{equation} \phi(X,Y,Z)= \phi(-X,Y,Z)=\phi(X,-Y,Z)=-\phi(X,Y,-Z). \end{equation} In particular, the fact that $\phi$ is odd in $Z$ means that the matching condition reduces to the condition that $\phi_{Z}$ vanishes on $W_{\infty}$. We show this by an indirect argument. With 
$$F(w,u)= Q^{-3/2} (Q+u^{2}) \tan^{-1}(u/\sqrt{Q}) + \kappa_{\epsilon}\  u $$ we have
$$ \frac{\partial F}{\partial u}= 2 u Q^{-3/2} \tan^{-1}(u/\sqrt{Q}) + (Q^{-1} + \kappa_{\epsilon}). $$ If $u$ is regarded as a function of $w,Z,A$ we have $\frac{\partial u}{\partial Z}=1$ so we get
\begin{equation}  \phi_{Z}= (2\pi i)^{-1} \int_{S^{1}} 2 Q^{-3/2} u \tan^{-1}(u/\sqrt{Q}) dw/w + K \end{equation}
where $K$ is a constant, independent of $(X,Y,Z)$.

Similarly, differentiating twice more, we get
\begin{equation}  \phi_{ZZZ}= (2/\pi i ) \int_{S^{1}} \frac{1}{(Q+ u^{2})^{2} } \  dw/w . \end{equation}
Since $\frac{\partial u}{\partial A}= w$ and $\frac{\partial u}{\partial \overline{A}}= w^{-1}$ we see that $\phi_{ZA\oA}=\phi_{ZZZ}$ (which is the Laplace equation for $\phi_{Z}$) and 
\begin{equation} \phi_{ZAA}= (2/\pi i) \int_{S^{1}} \frac{1}{(Q+ u^{2})^{2}}  w dw \ \ ,\ \  ,\ \ \  \phi_{Z\oA\oA}= (2/\pi i) \int_{S^{1}} \frac{1}{(Q+ u^{2})^{2} } w^{-3} dw .  \end{equation}

These equations are valid when the point $p=(X,Y,Z)$ does not lie in $\overline{W_{\infty}}$

\begin{lem}
The functions $\phi_{ZAA},\phi_{ZA\oA}, \phi_{Z\oA \oA}$ tend to $0$ as $p$ tends to a point of $W_{\infty}$. 
\end{lem}

To see this we need to consider the roots of the Laurent polynomial in $w$ $$ Q+u^{2}= 1+\epsilon/2(w^{2}+ w^{-2})+ (A w + Z -\oA w^{-1})^{2}, $$
when $Z$ is close to $0$. These are the roots of a quartic polynomial in $w$, so could be written explicitly in terms of radicals but we do not need that. Notice that the roots come in antipodal pairs. Write $\sigma$ for the antiholomorphic map $\sigma(w)= -\overline{w}^{-1}$: then if $w$ is a root so also is $\sigma(w)$.

 When $Z=0$ the polynomial can be written in terms of $v=w^{2}$:
    $$Q+u^{2}= (A^{2}+ \epsilon/2)v + (1-2\vert A\vert^{2}) + (\oA^{2}+\epsilon/2) v^{-1}. $$
There is an exceptional case when $A^{2}+\epsilon/2=0$. Then $1-2\vert A\vert^{2}=1-\vert \epsilon\vert >0$ so there are no roots. Otherwise, the roots are given by the quadratic formula
$$  (2A^{2}+\epsilon)^{-1} \left( (2\vert A\vert^{2}-1)\pm \sqrt{\Delta}\right) $$
where $$\Delta = (2\vert A\vert^{2}-1)^{2}-\vert 2A^{2}+\epsilon\vert^{2}= (1-\epsilon^{2})- 4 (\vert A\vert^{2}+ \epsilon {\rm Re}(A^{2}). $$
The equation $\Delta=0$ defines the ellipse $\Gamma$. For $p\in W_{0}$ we have $\Delta>0$ and the roots are real multiples of $(2A^{2}+\epsilon)^{-1}$. The exceptional case when $A^{2}+\epsilon/2=0$ is a limiting situation when the roots move to $0$ and $\infty$. But the relevant case for us is when $p\in W_{\infty}$, so $\Delta$ is negative and there are two distinct  roots which lie on the unit circle, say $v_{1}, v_{2}$. Thus for $p\in W_{\infty}$ the four roots of $Q+u^{2}$ can be written
$w_{1},-w_{1}, w_{2}, -w_{2}$, where $w_{1},w_{2}$ lie on the unit circle and $w_{2}\neq \pm w_{1}$. Fix $p_{0}\in W_{\infty}$ and consider $p\in \bR^{3}$ close to $p_{0}$. Then the four roots of $Q+u^{2}$ can be written as
$$ w_{1}(p)\ , \ \sigma(w_{1}(p))\ ,\  w_{2}(p)\ ,\  \sigma(w_{2}(p)) $$
where $w_{1}(p), w_{2}(p)$ vary continously with $p$ and tend to $w_{1}, w_{2}$ as $p\rightarrow p_{0}$. We choose $w_{1}(p), w_{2}(p)$ to lie {\it outside} the unit circle (for $p$ not in $W_{\infty}$) so $\sigma(w_{1}), \sigma(w_{2})$ lie {\it inside} the circle.

For $k=-3,-1,1$  let $\Omega_{k}$ be the meromorphic $1$-form
$$  \Omega_{k}= \frac{1}{(Q+u^{2})^{2}} w^{k} dw, $$
depending on $p$ close to $p_{0}$ as above. 
 This has poles  at the four roots of $Q+u^{2}$ and one checks that it has no poles at $0,\infty$. To simplify notation, fix attention first on the case $k=-1$. Write ${\rm Res}\  w_{1}(p)$ etc. for the residues  of $\Omega_{-1}$ at the four roots; these clearly depend continuously on $p$. By the residue theorem we have, for $p$ close to $p_{0}$ but not in $W_{\infty}$:
\begin{equation}  \phi_{ZA\oA}= (2/\pi i) \int_{S^{1} } \Omega_{-1} = -4 i ( {\rm Res}\ w_{1}(p)+ {\rm Res}\ w_{2}(p)= +4i ( {\rm Res}\ \sigma(w_{1}(p))+ {\rm
Res}\ \sigma(w_{2}(p)). \end{equation}
So, taking the limit as $p\rightarrow p_{0}$, we want to show that
\begin{equation} {\rm Res}\  w_{1} (p_{0})+ {\rm Res}\ w_{2}(p_{0})=0 \end{equation}
We have $w_{1}(p_{0})=w_{1}$ and  $\sigma(w_{1})=-w_{1}$. When $p=p_{0}$ the function $Q+u^{2}$ can be written in terms of $v=w^{2}$ as above and
$w^{-1}dw= \frac{1}{2} v^{-1} dv$. It follows that the residues of $\Omega_{-1}$ at $w_{1}$ and $-w_{1}$ are equal and similarly for $\pm w_{2}$. But taking the limit as $p\rightarrow p_{0}$ in (12)  we have
$$    {\rm Res}\ w_{1}+
{\rm Res} \ w_{2}= - \left[ {\rm Res}\ ( -w_{1})+ {\rm
Res}\ (-w_{2}) \ \right], $$
so both sides must vanish, as required. The argument for the other derivatives, using the $1$-forms $\Omega_{1},\Omega_{-3}$ is exactly the same.

\

Write $\psi$ for the restriction of $\phi_{Z}$ to $W_{\infty}$. We have shown that all the second derivatives of $\psi$ vanish, so $\psi$ is an affine linear function $\psi(X,Y)= h_{1} X + h_{2} Y + h_{3}$. But the symmetry (8) of $\phi$ implies that $\psi$ is even in $X,Y$,  so $h_{1}, h_{2}$ vanish and $\psi$ is a constant. By choosing the constant $\kappa_{\epsilon}$ appropriately we can arrange that $\psi=0$. Then $\phi$ satisfies the matching condition required to define a multivalued function,  as claimed.

\

 The next two Propositions complete the proof of Theorem 1. 
 
\begin{lem}
The derivatives $\phi_{Z}, \phi_{A}$ are bounded in a neighbourhood of $\Gamma$.
\end{lem}
We have seen that $$\phi_{Z}=(\pi i)^{-1} \int_{S^{1} }  Q^{-3/2} u \tan^{-1}(u\sqrt{Q}) dw/w + K, $$  and similarly one gets
$$  \phi_{A}= (\pi i)^{-1} \int_{S^{1} }   Q^{-3/2} u \tan^{-1}(u\sqrt{Q}) dw. $$
On the unit circle the function $Q$ is real and bounded above and below
   $$     1+\epsilon\geq Q \geq  1-\epsilon>0. $$
  For points $p$ in a neighbourhood of $\Gamma$ the function $u$ is bounded on the unit circle so we get
$$  \vert \phi_{A}\vert \leq C_{1} \int_{S^{1}} \vert \tan^{-1}(u/\sqrt{Q}) \vert $$ for some $C_{1}$ depending only on $\epsilon$,  and similarly for $\phi_{Z}$. 

Recall that
$$  \tan^{-1}(\zeta) = \frac{1}{2i} \log\left(\frac{i-\zeta}{i+\zeta}\right) $$
so $$ \vert \tan^{-1}(u/\sqrt{Q}) \leq \frac{1}{2} \left( \vert \log (\sqrt{Q}- iu)  \vert + \vert \log(\sqrt{Q} + iu)\vert\right). $$
On the other hand
$$  \log (Q+u^{2})= \log(\sqrt{Q}- iu) + \log(\sqrt{Q}+i u). $$
On the unit circle at most one of $\sqrt{Q}\pm iu$ is small in modulus so it follows that we have a bound
$$  \vert \tan^{-1}(u/\sqrt{Q})\vert \leq C_{2} \vert \log(Q+u^{2})\vert + C_{3}, $$
for constants $C_{2}, C_{3}$. Now the Lemma follows easily from the fact that
 that for any Laurent polynomial $P$ the integral of $\vert \log P\vert$ around the unit circle is finite, with a bound that is uniform over any compact set of such polynomials. 
\begin{prop}
The function $\tphi={\rm Sign} (Z)\ \phi$ satisfies the asymptotic condition
(2) of Theorem 1.
\end{prop}
We work with the original notation
\begin{equation}\phi= (2\pi)^{-1} \int_{0}^{2\pi} Q^{-3/2} \tan^{-1}(u/\sqrt{Q}) (Q+u^{2}) d\theta + \kappa_{\epsilon} Z,   \end{equation}
where $u = Z + i(X\cos\theta + Y\sin \theta)$,  and consider points with $Z>0$, so $\phi=\tphi$. We write $R=\sqrt{X^{2}+Y^{2}+Z^{2}}$. We compare with the simpler integral
\begin{equation}\uphi= (1/4)  \int_{0}^{2\pi}  Q^{-3/2}  (Q+u^{2}) d\theta.
\end{equation}
For ${\rm Re}\ \zeta>0$ and $\vert \zeta \vert \geq 1$ the function $\tan^{-1}$ satisfies an estimate
$$  \vert \tan^{-1}(\zeta)- \pi/2\vert \leq C_{4} \vert \zeta\vert^{-1} $$
and, for all $\zeta$, 
   $$  \vert (1+ \zeta^{2} \tan^{-1}(\zeta)\vert \leq C_{5} (1+\vert \zeta\vert^{2}). $$
Recall that $Q$ and $Q^{-1}$ are bounded on the unit circle. A moments thought shows that there are constants $C_{6}, C_{7}>0$ such that 
$\vert u/\sqrt{Q}\vert\geq C_{6} R$ except on a set $I\subset S^{1}$ of measure at most $C_{7}/R$. It follows that the contributions from $I$ to the integrals in (14),(15) are $O(R)$ as also is the term $\kappa_{\epsilon} Z$.
On the other hand
$$  \int_{S^{1}\setminus I} Q^{-3/2} (Q+u^{2}) \left(\tan^{-1}(u/\sqrt{Q})- \pi/2\right)  d\theta, $$
is bounded in modulus by 
$$  C_{8} R^{-1} \int_{S^{1}\setminus I}  Q^{-3/2}\vert Q+u^{2}\vert d\theta, $$which is also $O(R)$. So $\phi$ and $\uphi$ agree to $O(R)$, for large $R$. 
Now it is clear that $\uphi$ is asymptotic to 
$$  (1/4)\int_{0}^{2\pi} u^{2} Q^{-3/2} \ d\theta=-\lambda X^{2} - \mu Y^{2} + \nu Z^{2}. $$

\section{Further discussion and results}
\subsection{Existence and uniqueness}
The case when the parameter $\epsilon$ is zero is special because the solution is rotationally invariant. In cylindrical co-ordinates $(r,Z)$ with $r^{2}=X^{2}+Y^{2}$ we have
$$\phi= (2\pi)^{-1} \int_{0}^{2\pi} f(Z+i r \cos \theta) d\theta, $$
where $f(\zeta)= (1+\zeta^{2}) \tan^{-1}\zeta + \kappa_{\epsilon} \zeta$.

The three functions $\lambda(\epsilon),\mu(\epsilon), \nu(\epsilon)$ of Theorem 1 satisfy
$$  \lambda(\epsilon) +\mu(\epsilon) =\nu(\epsilon) \ \ \ \ \nu(-\epsilon)=\nu(\epsilon) \ \ \  \lambda(-\epsilon)=\mu(\epsilon). $$
Set $\varpi(\epsilon)= \lambda(\epsilon)/\nu(\epsilon)$ so $\varpi(-\epsilon)=1-\varpi(\epsilon)$ As $\epsilon$ tends to $-1$ the integrand defining $\nu(\epsilon)$ converges pointwise to $2^{-3/2} \vert \sin \theta \vert^{-3}$ while that defining $\mu(\epsilon)$ converges to $2^{-3/2} \vert \sin \theta \vert^{-1}$. Straightforward arguments show that $\nu(\epsilon)$ is $O((1-\epsilon)^{-2})$ as $\epsilon\rightarrow 1$ whereas $\mu(\epsilon)$ is $O(\vert \log(1-\epsilon)\vert$. Thus $\varpi(\epsilon)$ tends to $1$ as $\epsilon\rightarrow -1$ and to $0$ as $\epsilon\rightarrow 1$. 

\begin{lem}

The function $\varpi$ is strictly decreasing on the interval $(-1,1)$.
\end{lem}

We consider a slightly more general situation. Suppose that $\alpha>0$ and  $f$ is a positive function, not equal to a constant.  For  $\eta>0$ define
$$I(\eta)= \int \frac{1}{(1+\eta f)^{\alpha}} d\theta , J(\eta)= \int \frac{f}{(1+\eta f)^{\alpha}}  d\theta. $$
We claim that the ratio $J/I$ is a decreasing function of $\eta$. Let $dm$ be the measure 
$$   dm= \frac{1}{(1+\eta f)^{1+\alpha}} \ d\theta. $$
Elementary calculation shows that
$$  J \frac{dI}{d\eta} - I \frac{dJ}{d\eta} = \alpha \left[ \int f^{2} dm \int dm - \left( \int f dm \right)^{2}\right]. $$
The right hand side is strictly positive by the Cauchy-Schwartz inequality and the left hand side is $-I^{2}$ times the derivative of $J/I$.

In the case at hand take $f= \cos ^{2}\theta$ and (for $0<\epsilon<1$) set $\eta= 2\epsilon/(1-\epsilon)$. Then
$\varpi(\epsilon)= J(\eta)/I(\eta)$ so $\varpi(\epsilon)$ is decreasing over this range. The symmetry $\varpi(-\epsilon)=1-\varpi(\epsilon)$ handles the range of negative $\epsilon$.  

\

\

Our solutions for parameter values $\pm \epsilon$  are geometrically equivalent: they correspond under interchange of the $X,Y$ co-ordinates. We can obviously obtain more multivalued functions by applying dilations and Euclidean motions and multiplying by constants. If $\phi$ is asymptotic to a quadratic function $Q$ at infinity then the rescalings
\begin{equation}   \rho^{-2}\phi( \rho p), \end{equation}
for non-zero $\rho$ have the same asymptotics, as also do the translations.    It follows from Lemma 4 that  any non-degenerate quadratic form on $\bR^{3}$ with trace zero is equivalent under scaling and orthogonal transformations to $\lambda(\epsilon) X^{2}+ \mu(\epsilon) Y^{2}- nu(\epsilon) Z^{2}$ for a unique $\epsilon>0$. Thus we have:
\begin{prop}
For any trace-free, nondegenerate quadratic form $Q$ on $\bR^{3}$ there is a multivalued harmonic function with $\tphi$ asymptotic to $Q$. Moreover within the family of solutions we have constructed, the solution is unique up to translation and the rescaling (16).
\end{prop}

It is an interesting question whether the solution, with given quadratic asymptotics, is absolutely unique: in other words whether there is some quite different family of solutions. In the case of circular symmetry this is proved by Yan in \cite{kn:Yan}.

\subsection{Complex geometry}

In the twistor setting the key to our construction is  the complex curve
$\Sigma_{\epsilon} \subset \bT$ defined in our co-ordinate patch by the equation $Q(w)+u^{2}=0$ (extended to the whole of $\bT$ in the standard way). For $\epsilon\neq
0$ it is a smooth curve of genus $1$ and when $\epsilon=0$ it degenerates
into two sections $u=\pm i$ of $\bT\rightarrow S^{2}$, corresponding to two
purely imaginary points in $\bC^{3}$. In this subsection we discuss the relation between $\Gamma\subset \bR^{3}$, which is the singular set of  the multivalued function, and
$\Sigma_{\epsilon}$ which is the singular set of the cohomological data on $\bT$. 

Recall that $\bC^{3}$ is identified with the set of sections of $\bT\rightarrow S^{2}$. In the complexified picture there is an alternative definition of $\bT$ as a set of complex hyperplanes in $  \bC^{3}$. For the temporary purposes of this discussion,  say that a hyperplane $H$ is \lq\lq special'' if  the restriction of the quadratic differential $dX^{2}+dY^{2}+dZ^{2}$ to $H$ has rank $1$. Then $\bT$ can be viewed as the set of special hyperplanes. To relate this to the previous description, given an oriented line in $\bR^{3}$ choose an oriented orthonormal frame $\lambda,\mu$ for the normal bundle to the line. Then the line can be defined by  equations for $p\in \bR^{3}$:
$$  \lambda.p= c_{1}  \ \ \ \ \mu.p=c_{2}, $$
in terms of the usual dot product and for suitable $c_{1}, c_{2}$. The special hyperplane in $\bC^{3}$ corresponding to the line is the set of complex solutions $p\in \bC^{3}$ of the equation
$$   (\lambda+ i\mu). p = (c_{1} + i c_{2}). $$

For any curve $\Sigma\subset \bT$ we define a set $\Gamma^{\bC}\subset \bC^{3}$ to be the sections which are tangent to $\Sigma$ at some point. It is a complex hypersurface in $\bC^{3}$ and has the property that the tangent space at each point is special, in the sense above.   There is a dual relation between
$\Sigma$ and $\Gamma^{\bC}$.  Given any hypersurface $V\subset \bC^{3}$ with the property that all tangent spaces are special we can define a set $S[V]\subset \bT$ of special hyperplanes tangent at some point to $V$. Ignoring problems of singularities, $S[V]$ is the image of a \lq\lq Gauss map'' $G:V\rightarrow \bT$ and is a complex curve in $\bT$. (The one-dimensional fibres of $G$ foliate $V$  according to the null directions in the tangent spaces.) The duality relation is
$$   \Sigma= S[ \Gamma^{\bC}(\Sigma)]. $$

So, ignoring technicalities, there is a 1-1 correspondence between \lq\lq special' hypersurfaces in $\bC^{3}$ and curves in $\bT$. 

Going back to our specific case with the curves $\Sigma_{\epsilon}$, the ellipse $\Gamma_{\epsilon}$ is a connected component of the real points
$$  \Gamma^{\bC}(\Sigma_{\epsilon})\cap \bR^{3}. $$
The situation is a little more complicated than that sketched above, because $\Gamma^{\bC}$ is singular at these points: it has a normal crossing singularity corresponding to the fact that for $p\in \Gamma$ the section $L_{p}$ is tangent to $\Sigma_{\epsilon}$ at {\it two} points, interchanged by the real structure. The tangents to these two components define two maps
$$  G^{\pm}: \Gamma\rightarrow \Sigma_{\epsilon}. $$
Unwinding the definitions, these are the maps defined by the tangent lines to $\Gamma$,  for the two choices of orientation. The images correspond to the solutions of the equation $u^{2} + Q(w)=0$ for $w$ on the unit circle. The duality relation comes down to following   high school geometry exercise.
Let $L$ be a tangent line to the ellipse 
$$ (1+\epsilon)^{-1} X^{2} + (1-\epsilon)^{-1}Y^{2}=1 $$ which makes angle $\theta$ with the $Y$-axis. Let $P$ be the distance from $L$ to the origin then
$$   P^{2}= 1+\epsilon \cos 2\theta . $$ 

\subsection{Evaluating some integrals}

The functions $\lambda(\epsilon), \mu(\epsilon), \nu(\epsilon)$ are defined by periods of meromorphic differentials on the family of elliptic curves $\Sigma_{\epsilon}$. They satisfy second order differential equations with regular singularities at the three points $\epsilon= \pm 1, \infty$. These differential equations can be transformed to hypergeometric equations and the functions $\lambda,\mu,\nu$ expressed in terms of standard hypergeometric functions.  

\

In the case when $\epsilon=0$ we can pin down the constant $\kappa=\kappa_{0}$ in (7). In that case we have
$$  \phi_{Z}= (2\pi)^{-1}\int_{0}^{2\pi} f'(Z+i r \cos \theta) d\theta +\kappa $$
where $f(\zeta)= (1+\zeta^{2}) \tan^{-1}(\zeta)$ so 
$$  f'(\zeta)= 1+ 2\zeta \tan^{-1}(\zeta $$
We take $Z=0$ and $r<1$ (corresponding to a point of $W_{0}$). Then
\begin{equation}\phi_{Z}= 1 + \kappa + \pi^{-1} \int_{0}^{2\pi} i r\cos\theta \tan^{-1}(ir\cos \theta) d\theta. \end{equation}

Change variable to $t= r \cos \theta$; then the last term on the right hand side of (17) is
$$ I= (2i/\pi) \int_{-r}^{r} \frac{t}{\sqrt{r^{2}-t^{2}}}  \tan^{-1}(it) dt . $$
We integrate by parts to write this as
$$I= -(2/\pi) \int_{-r}^{r} \sqrt{r^{2}-t^{2}} \frac{1}{1-t^{2}} dt. $$
Taking $t$ to be a complex variable this is
\begin{equation} I= -\pi^{-1} \int_{C} \frac{\sqrt{r^{2}-t^{2}}}{1-t^{2}}  dt, \end{equation}
where the square root is defined by taking a cut on the interval $[-r,r]$ and $C$ is a closed contour around this interval, not enclosing the points $t=\pm 1$. Keeping track of the sign in the square root one sees that the integrand in (18) is asymptotic to $i t^{-1}$ at infinity, so the integral around a large circle is $-2\pi$. Calculating the residues at the poles $t-\pm 1$ gives  $I=  4\pi(\sqrt{1-r^{2}}-1)$ so, finally,
$$   \phi_{Z}= \kappa-1 + 2 \sqrt{1-r^{2}}. $$
Taking the limit as $r\rightarrow 1$ we see that we need $\kappa=1$ to make $\phi_{Z}$ vanish on $W_{\infty}$. The author does not know any simple expression for $\kappa$ when $\epsilon\neq 0$.

%+Bibliography

%-Bibliography

\end{document}